\documentclass{jnmp}

\usepackage{amsmath}

\JNMPnumberwithin{equation}{section}

\newtheorem{theorem}{Theorem}
\newtheorem{lemma}{Lemma}

\theoremstyle{definition}
\newtheorem{definition}{Definition}
\newtheorem{remark}{Remark}

\begin{document}

\renewcommand{\evenhead}{S~Evje and K~H~Karlsen}
\renewcommand{\oddhead}{An Error Estimate for Viscous Approximations}

\thispagestyle{empty}

\FirstPageHead{9}{3}{2002}{\pageref{evije-firstpage}--\pageref{evije-lastpage}}{Article}

\copyrightnote{2002}{S~Evje and K~H~Karlsen}

\Name{An Error Estimate for Viscous
Approximate  Solutions of
Degenerate Parabolic Equations}
\label{evije-firstpage}

\Author{Steinar EVJE~$^{\dag}$ and Kenneth H~KARLSEN~$^{\ddag}$}

\Address{$^\dag$~RF-Rogaland Research, Thorm{\o}hlensgt.~55, N--5008 Bergen, Norway\\
~~E-mail: Steinar.Evje@rf.no\\[10pt]
$^{\ddag}$~Department of Mathematics, University of Bergen,\\
~~Johs.\ Brunsgt.\ 12, N--5008 Bergen, Norway\\
~~E-mail: kennethk@math.uib.no, \ \ URL: http://www.mi.uib.no/\~{}kennethk/}

\Date{Received February 14, 2001; Revised  December 13, 2001;
Accepted February 13, 2002}

\begin{abstract}
\noindent
  Relying on  recent advances in the theory of entropy
   solutions for nonlinear (strongly) degenerate parabolic
   equations, we present a \textit{direct} proof
   of an $L^1$ error estimate for viscous approximate
   solutions of the initial value problem for
   \[\partial_t w+\mathrm{div} \bigl(V(x)f(w)\bigr)  = \Delta A(w),\] where
   $V=V(x)$ is a vector field, $f=f(u)$ is a scalar function,
   and $A'(\cdot)\ge 0$. The viscous approximate
   solutions are weak solutions of the initial value problem for
   the uniformly parabolic equation
   \[\partial_t w^{\varepsilon}+\mathrm{div} \bigl(V(x)f(w^{\varepsilon})\bigr)  =
   \Delta \bigl(A(w^{\varepsilon}) +\varepsilon w^{\varepsilon}\bigr), 
   \qquad \varepsilon>0.
   \]
   The error estimate is of order $\sqrt{\varepsilon}$.
\end{abstract}

\section{Introduction}
In this paper we are interested in certain ``viscous'' approximations of
entropy solutions of the initial value problem
\begin{gather}
      \partial_t w+\mathrm{div} \bigl(V(x)f(w)\bigr)  = \Delta A(w), \qquad (x,t)\in Q_T, \nonumber\\
      w(x,0)=w_{0}(x), \qquad  x\in {\mathbb R}^d,   \label{cont_model}
\end{gather}
where $Q_T={\mathbb R}^d\times (0,T)$ with $T>0$ fixed, $u:Q_T\to {\mathbb R}$ is
the sough function, $V:{\mathbb R}^d\to{\mathbb R}$ is a
(not necessarily divergence free) velocity
field, $f:{\mathbb R}\to{\mathbb R}$ is the convective flux
function, and $A:{\mathbb R}\to{\mathbb R}$ is the ``diffusion'' function.
For the diffusion function the basic assumption is that
$A(\cdot)$ is nonincreasing.
This condition implies that \eqref{cont_model} is
a (strongly) \textit{degenerate parabolic} problem. For example,
the hyperbolic equation
$\partial_t w+\mathrm{div} \bigl(V(x)f(w)\bigr) =0$ is a~special case of~\eqref{cont_model}.
Problems such as~\eqref{cont_model} occur in several
important applications. We mention here only
two examples: flow in porous
media (see, e.g., \cite{EspKar}) and
sedimentation-consolidation processes~\cite{Burger:bok}.

Since $A(\cdot)$ is merely nondecreasing, solutions are not
necessarily smooth and weak solutions must be sought. Moreover, as
is well-known in the theory of hyperbolic conservation laws, weak
solutions are not uniquely determined by their initial data. To
have a~well-posed problem we need to consider entropy solutions,
i.e., weak solutions that satisfy a~Kru\v{z}kov--Vol'pert type
entropy condition. A precise statement is given in
Section~\ref{sec:result} (see Definition~\ref{def:sol}). For
purely hyperbolic equations this entropy condition was introduced
by Kru\v{z}kov~\cite{Kruzkov} and Vol'pert~\cite{Volpert}. For
degenerate parabolic equations, it was introduced by Vol'pert and
Hudjaev~\cite{VolHud}.

Following Carrillo~\cite{Carrillo}, Karlsen and Risebro~\cite{KR}
proved that the entropy solution of~\eqref{cont_model} (as well as
a more general equation) is unique. Moreover in the
$L^{\infty}\big(0,T;BV\big({\mathbb R}^d\big)\big)$ class of
entropy solutions, they proved an $L^1$ contraction principle.
Existence of an $L^{\infty}\big(0,T;BV\big({\mathbb
R}^d\big)\big)$ entropy solution of~\eqref{cont_model} follows
from the results in Vol'pert and Hu\-djaev~\cite{VolHud} or
Karlsen and Risebro~\cite{KR:Rough_Diff} (the latter deals with
convergence of finite difference methods). The proof in~\cite{KR}
of uniqueness and stability is based on the ``doubling of
variables'' strategy introduced in Carrillo~\cite{Carrillo} (see
also Chen and DiBenedetto~\cite{ChenDiBen}), which in turn is a
generalization of the pioneering work by
Kru\v{z}kov~\cite{Kruzkov} on hyperbolic equations. Related papers
dealing with the ``doubling of variables'' device for degenerate
parabolic equations include, among others,
Carrillo~\cite{Carrillo:94}, Otto~\cite{Otto:L1_Contr}, Rouvre and
Gagneux~\cite{RouvreGagneux}, Cockburn and
Gripenberg~\cite{CockGripen} B\"urger, Evje and
Karlsen~\cite{BurgerEvjeKarlsen:1D_IBVP,BK:Hyp2000},
Ohlberger~\cite{Ohlberger:FVM}, Mascia, Porretta, and
Terracina~\cite{Mascia_etal:2000}, Eymard, Gallouet, Herbin and
Michel~\cite{Eymard_etal_II:2000}, and Karlsen and
Ohlberger~\cite{KO}.

In this paper we are interested in certain
approximate solutions of~\eqref{cont_model} coming from
solving the uniformly parabolic problem
\begin{gather}
        \partial_t w^{\varepsilon}+\mathrm{div} \bigl(V(x)f(w^{\varepsilon})\bigr)
      = \Delta A^{\varepsilon}(w^{\varepsilon}), \qquad (x,t)\in Q_T,\nonumber\\
      w^{\varepsilon}(x,0)=w_{0}(x), \qquad x\in {\mathbb R}^d,
\label{viscous_model}
\end{gather}
where $A^{\varepsilon}(w^{\varepsilon})=A(w^{\varepsilon})+\varepsilon w^{\varepsilon}$,
 $\varepsilon>0$.
We refer to $w^{\varepsilon}$ as a \textit{viscous} approximate
solution of~\eqref{cont_model}.
Convergence of $w^{\varepsilon}$ to the unique entropy
solution $w$ of \eqref{cont_model} as $\varepsilon\downarrow 0$ follows
from the results in Vol'pert and
Hudjaev~\cite{VolHud}.
Our main interest here is to give an explicit rate of convergence
for $w^{\varepsilon}$ as $\varepsilon\downarrow 0$, i.e., an $L^1$ error estimate for
viscous approximate solutions.

There are several ways to prove such an error estimate.
One way is to view it as a consequence of a continuous dependence estimate.
Combining the ideas in~\cite{KR} with those in Cockburn
and Gripenberg~\cite{CockGripen}, who used
a variant of Kru\v{z}kov's ``doubling of variables'' device for~\eqref{cont_model}
with $V\equiv 1$, Evje, Karlsen and Risebro~\cite{EKR} established an explicit
``continuous dependence on the nonlinearities'' estimate
for entropy solutions of~\eqref{cont_model}.
A~direct consequence of this estimate is the
error bound $\|w^{\varepsilon}-w\|_{L^1|(Q_T)}
=\mathcal{O}\bigl(\sqrt{\varepsilon}\bigr)$, at least
when  $w^{\varepsilon},w$ belong to $L^{\infty}\left(0,T;BV\left({\mathbb R}^d\right)\right)$ and
$V$ is sufficiently regular.
Unfortunately the techniques employed in~\cite{EKR} require that
one works with (smooth) viscous approximations of~\eqref{cont_model}.
The proof in~\cite{EKR} (as well as the one in~\cite{CockGripen})
did not exploit the entropy solution ``machinery''
developed by Carrillo~\cite{Carrillo}.

The main purpose of this work is to show that one can indeed use
the ``doubling of variables'' device to compare \textit{directly} the
entropy solution~$w$ of~\eqref{cont_model} against the viscous
approximation $w^{\varepsilon}$ of~\eqref{viscous_model}.
Hence there is no need to work with
approximate solutions of~\eqref{cont_model}.
Although our proof is of independent
interest, it may also shed some light on
how to obtain error estimates for numerical methods. Most
numerical methods (related to this class of equations) have~\eqref{viscous_model} as
a ``model'' problem and, in this context, the size of~$\varepsilon$ designates
the amount of ``diffusion'' present in the numerical method.
A step in the direction of obtaining error
estimates for numerical methods has been
taken by Ohlberger~\cite{Ohlberger:FVM} with his a posteriori
error estimate for a finite volume method.
We will in future work use the ideas devised herein
to derive error estimates a priori for finite difference methods.

The rest of this paper is  organized as follows: In Section \ref{sec:result}
we state the definition of an entropy solution
and the main result (Theorem~\ref{main_thm}).
Section~\ref{sec:entropy} is
devoted to the derivation of certain entropy inequalities
for the exact entropy solution and its viscous
approximation. Equipped with these
entropy inequalities, we prove
the error estimate (Theorem~\ref{main_thm})
in Section~\ref{sec:main_thm}.

\section{Statement of result}\label{sec:result}

Following~\cite{KR:Rough_Diff,KR} we start by stating sufficient
conditions on $V=(V_1,\dots,V_d),f,A$ and $u_0$ to ensure the
existence of a unique $L^{\infty}(0,T;BV({\mathbb R}^d))$
entropy solution of~\eqref{cont_model}:
\begin{gather}
         V\in \big(L^{\infty}\big({\mathbb R}^d\big)\big)^d\cap
      \big(\mathrm{Lip}\big({\mathbb R}^d\big)\big)^d;\quad
      \mathrm{div} V\in BV\big({\mathbb R}^d\big);\nonumber\\
      f\in \mathrm{Lip_{loc}}({\mathbb R});\quad f(0)=0;\nonumber\\
      A\in \mathrm{Lip_{loc}}({\mathbb R})\quad
      \text{and}\quad \text{$A(\cdot)$ is nondecreasing with $A(0)=0$};\nonumber\\
      u_0\in L^{\infty}\big({\mathbb R}^d\big)\cap BV\big({\mathbb R}^d\big).
 \label{Coeff_cond}
\end{gather}
Note that the first
condition in~\eqref{Coeff_cond} implies
\[
V\in \big(W^{1,1}_{\mathrm{loc}}\big({\mathbb R}^d\big)\big)^d.
\]

In \eqref{Coeff_cond} and elsewhere in
this paper the space $BV\big({\mathbb R}^d\big)$ is defined as
\[
BV\big({\mathbb R}^d\big)=\left\{g\in L^1\big({\mathbb R}^d\big):
 |g|_{BV\left({\mathbb R}^d\right)}<\infty\right\},
\]
where $|g|_{BV\big({\mathbb R}^d\big)}$ denotes the total variation
of $g$, i.e., $g\in BV\big({\mathbb R}^d\big)$ if and only if
$g\in L^1\big({\mathbb R}^d\big)$ and the first
order distributional derivatives of $g$ are represented
by finite measures on ${\mathbb R}^d$.

Equipped with \eqref{Coeff_cond} we can state
the following definition of an entropy solution:
\begin{definition}[Entropy Solution]
   \label{def:sol}
   A function $w(x,t)$ is called an entropy solution \linebreak of~\eqref{cont_model} if
   \begin{itemize}
      \item[(i)]
        $w\in L^1(Q_T) \cap L^{\infty}(Q_T)\cap
        C\left(0,T;L^1\left({\mathbb R}^d\right)\right)$,
      \item[(ii)]
      $A(w)\in L^2\left(0,T;H^1\left({\mathbb R}^d\right)\right)$,
      \item[(iii)]
      $w(x,t)$ satisfies the {\it entropy inequality}
      \begin{gather}
            \iint\limits_{Q_T}\Bigl(
            |w-k|\partial_t \phi+\mathrm{sgn}(w-k)\bigl[ V(x)(f(w)-f(k))-
            \nabla A(w) \bigr]\cdot\nabla\phi\nonumber\\
 \qquad {}- \mathrm{sgn}(w-k)\mathrm{div} V(x) f(k)\phi\Bigr)\, dt\, dx\geq 0,
            \qquad \forall \; k\in{\mathbb R},         \label{entropy_ineq}
      \end{gather}
      for all nonnegative $\phi\in C^{\infty}_0(Q_T)$ and
      \item[(iv)]
      $\|w(\cdot,t)-w_0\|_{L^1\left({\mathbb R}^d\right)} \rightarrow 0$
      as $t\downarrow 0$ (essentially).
   \end{itemize}
\end{definition}

Note that, if we take $k>{\rm ess}\sup w(x,t)$ and $k<{\rm ess}\inf
w(x,t)$ in \eqref{entropy_ineq}, then an
approximation argument reveals that
\begin{equation}
   \label{L2_weak_tmp}
   \iint\limits_{Q_T}\Bigl(w\phi_t +
   \bigl[V(x)f(w) - \nabla A(w)\bigr]\cdot\nabla\phi \Bigr)\, dt\, dx  = 0
\end{equation}
holds for all $\phi\in H^1(Q_T)$.
Let $\langle\cdot,\cdot\rangle$ denote the usual pairing
between $H^{-1}\big({\mathbb R}^d\big)$ and $H^1\big({\mathbb R}^d\big)$.
 From \eqref{L2_weak_tmp} we conclude  that
\[
\partial_t  w\in L^2\big(0,T;H^{-1}\big({\mathbb R}^d\big)\big),
\]
so that
\begin{gather}
   \label{L2_weak}
   -\! \int_0^T\!\! \big\langle \partial_t  w,\phi\bigr\rangle \, dt
   \! +\!\! \iint\limits_{Q_T}
   \Bigl(\bigl[V(x)f(w) - \nabla A(w)\bigr]
   \cdot\nabla\phi \Bigr)\, dt\, dx  = 0, \!\quad
   \forall \;\phi\in H^1(Q_T).\!\!
\end{gather}
In other words an entropy solution $w(x,t)$ of \eqref{cont_model} is also
a \textit{weak} solution of the same problem.

In this paper we are interested in comparing the
entropy solution~$w$ of~\eqref{cont_model} against the
weak solution~$w^{\varepsilon}$ of the viscous problem~\eqref{viscous_model}.
 From the results in Karlsen and Risebro~\cite{KR:Rough_Diff}
or Vol'pert and Hudjaev~\cite{VolHud} there
exists a weak solution $w^{\varepsilon}\in L^{\infty}\big(0,T;BV\big({\mathbb R}^d\big)\big)$
of~\eqref{viscous_model}. Since $A^{\varepsilon}(\cdot)$ is increasing,  the
uniqueness result in Karlsen
and Risebro~\cite{KR} (see also Remark~\ref{viktig_remark} herein) tells
us that this weak solution is in fact a unique solution. Moreover from
the energy estimate we conclude that $w^{\varepsilon} \in
L^2\big(0,T;H^1\big({\mathbb R}^d\big)\big)$. Of course, if $V$, $f$, $A$, $u_0$ are smooth
enough, one can prove that the weak solution~$w^{\varepsilon}$
of~\eqref{viscous_model} is actually
a~classical ($C^{2,1}$) solution.  See, e.g., Vol'pert and Hudjaev~\cite{VolHud}.
Here it will be sufficient to know
that $w^{\varepsilon}$ belongs to 
$L^2\big(0,T;H^1\big({\mathbb R}^d\big)\big)$ (not $C^{2,1}$).

We are now ready to state our  main theorem:
\begin{theorem}[Error Estimate]
   \label{main_thm}
   Suppose that the conditions in \eqref{Coeff_cond} hold.
   Let $w\in L^{\infty}\big(0,T;BV\big({\mathbb R}^d\big)\big)$ be
   the unique entropy solution of \eqref{cont_model}
   and let $w^{\varepsilon}\in L^2\big(0,T;H^1\big({\mathbb R}^d\big)\big)\cap
   L^{\infty}\big(0,T;BV\big({\mathbb R}^d\big)\big)$ be the unique weak solution
   of~\eqref{viscous_model}.
   Then there exists a constant~$C$, independent of~$\varepsilon$, such that
   \begin{equation}
      \label{Error_estimate}
      \|w^{\varepsilon}-w\|_{L^1(Q_T)}
       \le C \sqrt{\varepsilon}.
   \end{equation}
\end{theorem}

\section{Entropy inequalities}
\label{sec:entropy}
In Section \ref{sec:main_thm} we follow the uniqueness proof of
Carrillo~\cite{Carrillo} to obtain an estimate of
the difference between~$w^{\varepsilon}$ and~$w$.
To this end it will be necessary to derive
two entropy inequalities for the exact solution~$w$ and
two approximate entropy inequalities for the
viscous solution~$w^{\varepsilon}$.
The purpose of this section is to derive
these inequalities. (See Lemma~\ref{entropy_hypar}
and Lemma~\ref{viscous_hypar} below.)

Note that, differently from the pure
hyperbolic case~\cite{Kruzkov}, we need to operate with
one additional entropy inequality (actually an
equality for the exact solution~$w$) taking
into account the parabolic (dissipation) mechanism in
the equation. Hence we introduce a~set~$H$ corresponding
to the regions where $A(\cdot)$ is ``flat'' and~\eqref{cont_model}
behaves hyperbolically. More
precisely, let $A^{-1}:{\mathbb R}\to {\mathbb R}$ denote the unique
left-continuous function which
satisfies $A^{-1}(A(u))=u$ for all $u\in {\mathbb R}$.
Then we define
\[
H=\Bigl\{r\in{\mathbb R}\,:\,A^{-1}(\cdot)\mbox{ is discontinuous at $r$} \Bigr\}.
\]
Since $A(\cdot)$ is a monotonic function, $H$ is at most countable.
The dissipation mechanism in the equation
is effective only in the ($x,t$) region corresponding to
the complement of $H$.

To prove Lemma~\ref{entropy_hypar} and Lemma~\ref{viscous_hypar}
below we need the following ``weak'' chain rule:
\begin{lemma}
   \label{chain_rule}
   Let $u:Q_T\to{\mathbb R}$ be a measurable function satisfying
   the four conditions
   \begin{enumerate}
       \item[{\rm (1)}]  $u\in L^1(Q_T)\cap L^{\infty}(Q_T)
                   \cap C\big(0,T;L^1\big({\mathbb R}^d\big)\big)$,

       \item[{\rm (2)}]  $u(0,\cdot)=u_0\in L^1\big({\mathbb R}^d\big)
        \cap L^{\infty}\big({\mathbb R}^d\big)$,

       \item[{\rm (3)}]  $\partial_t  u \in 
       L^2\bigl(0,T; H^{-1}\big({\mathbb R}^d\big)\bigr)$ and

       \item[{\rm (4)}]  $A(u)\in L^2\bigl(0,T;H^1\big({\mathbb R}^d\big)\bigr)$.
   \end{enumerate}
   For every nonnegative and compactly supported
   $\phi\in C^{\infty}(Q_T)$ with
   $\phi|_{t=0}=\phi|_{t=T}=0$ we have
   \[
   - \int_0^T \Bigl\langle \partial_t  u,\psi\bigl( A(u)\bigr) \phi
   \Bigr\rangle \, dt
   = \iint\limits_{Q_T}\biggl(\int_k^u \psi(A(\xi))\, dxi\biggr)\phi_t\, dt\, dx,
   \qquad k\in{\mathbb R},
   \]
   where  $\psi:{\mathbb R}\to{\mathbb R}$ is a nondecreasing and
   Lipschitz continuous function.
\end{lemma}
The proof of Lemma~\ref{chain_rule} is very
similar to the proof of the ``weak chain'' rule in
Carrillo~\cite{Carrillo} and it is therefore
omitted. See instead~\cite{KR}.

The following lemma, which deals with entropy inequalities for the exact entropy
solution~$w$, is a direct consequence of the
very definition of an entropy solution.
\begin{lemma}
   \label{entropy_hypar}
   The unique entropy solution $w$ of \eqref{cont_model} satisfies:
   \begin{itemize}
      \item[{\rm (i)}]
      For all $k\in{\mathbb R}$ and all nonnegative $\phi\in C^{\infty}_0(Q_T)$ we have
      \begin{gather}
         \label{hyp_entw}
         E^{\mathrm{hyp}}(w,k,\phi)\geq 0,
      \end{gather}
      where
      \begin{gather}
            E^{\mathrm{hyp}}(w,k,\phi) :=
            \iint\limits_{Q_T}\Bigl(  |w -k|\partial_t  \phi
            + \mathrm{sgn}(w -k)\bigl[ V(x)(f(w)- f(k))
            \nonumber\\
\phantom{E^{\mathrm{hyp}}(w,k,\phi) :=}{}
 -\nabla A(w)\bigr] \cdot\nabla\phi
 - \mathrm{sgn}(w-k)\mathrm{div} V(x) f(k)\phi \Bigr) \, dt\, dx.
   \label{Ehyp}
      \end{gather}
      We refer to \eqref{hyp_entw} as
      a hyperbolic entropy inequality.
      \item[{\rm (ii)}]
      For all $k$ such that $A(k)\notin H$ and all
      nonnegative $\phi\in C^{\infty}_0(Q_T)$ we have
      \begin{gather}
         \label{par_entw}
         E^{\mathrm{par}}(w,k,\phi)= 0,
      \end{gather}
      where
      \begin{gather}
            E^{\mathrm{par}}(w,k,\phi) :=
            \iint\limits_{Q_T}\Bigl(  |w -k|\partial_t  \phi
            + \mathrm{sgn}(w-k)\bigl[V(x)(f(w)- f(k))
            \nonumber\\
            \phantom{E^{\mathrm{par}}(w,k,\phi) := }{}
-  \nabla A(w)\bigr] \cdot\nabla\phi
- \mathrm{sgn}(w-k)\mathrm{div} V(x) f(k)\phi \Bigr) \, dt\, dx
           \nonumber\\
\phantom{E^{\mathrm{par}}(w,k,\phi) := }{}
 - \lim_{\eta\downarrow 0} \iint\limits_{Q_T}
            \bigl|\nabla A(w) \bigr|^2\mathrm{sgn}_{\eta}'\bigl( A(w) - A(k)
            \bigr) \phi\, dt\, dx.         \label{Epar}
      \end{gather}
      In \eqref{Epar} (and elsewhere in this paper) $\mathrm{sgn}_{\eta}$
      is the approximate sign function defined by
      \begin{gather}
         \label{sgnappr}
         \mathrm{sgn}_\eta (\tau) :=
         \begin{cases}
            \mathrm{sgn}(\tau)  & \mbox{if  $|\tau| >   \eta$,} \\
            \tau/\eta & \mbox{if  $|\tau| \leq \eta$,}
         \end{cases}
         \qquad \eta>0.
      \end{gather}
     We refer to \eqref{par_entw} as
     a parabolic entropy inequality.
   \end{itemize}
\end{lemma}

\begin{proof}
The first inequality \eqref{hyp_entw} is nothing but the entropy condition
for the entropy solution~$w$.  So there is nothing to prove.
We turn to the proof of the second inequality~\eqref{par_entw},
which borrows a lot from Carrillo~\cite{Carrillo} (see also~\cite{KR}).
In what follows we always let~$k$ and~$\phi$ be as in
the lemma and the approximate sign function
$\mathrm{sgn}_{\eta}(\cdot )$ is always the one defined in~\eqref{sgnappr}.

Since $w$ satisfies \eqref{L2_weak} and
$\bigl[\mathrm{sgn}_{\eta}(A(w)-A(k))\phi\bigr]\in L^2\big(0,T; H^1\big({\mathbb R}^d\big)\big)$, we have
\begin{gather*}
    -\int_0^T \Bigl\langle \partial_t  w,
    \mathrm{sgn}_{\eta}(A(w)-A(k))\phi\Bigr\rangle \, dt
    \\  \qquad
    {}+ \iint\limits_{Q_T} \Bigl( \bigl[V(x)(f(w)-f(k)) - \nabla A(w) \bigr]
    \cdot\nabla\bigl[\mathrm{sgn}_{\eta}(A(w)-A(k))\phi\bigr]
   \\  \qquad
   {}- \mathrm{div} V(x) f(k)\bigl[\mathrm{sgn}_{\eta}(A(w)-A(k))\phi\bigr]\Bigr)\, dt\, dx =0.
\end{gather*}
Introduce the function $\psi_{\eta}(z) =
\mathrm{sgn}_{\eta}\bigl(z-A(k)\bigr)$ and note that Lemma \ref{chain_rule} can be
applied so that
\begin{gather*}
    -\int_0^T \Bigl\langle \partial_t  w,
    \mathrm{sgn}_{\eta}(A(w)-A(k))\phi\Bigr\rangle \, dt
    = \iint\limits_{Q_T} \biggl(\int_k^w \mathrm{sgn}_{\eta}(A(\xi)-A(k))\, d\xi
    \biggr) \partial_t \phi\, dt\, dx.
\end{gather*}
Hence
\begin{gather}
     \iint\limits_{Q_T} \biggl(\int_k^w \mathrm{sgn}_{\eta}(A(\xi)-A(k))\, d\xi
     \biggr)\partial_t  \phi \, dt\, dx
     \nonumber\\ \qquad
     {}+\iint\limits_{Q_T} \Bigl(\bigl[V(x)(f(w)-f(k)) - \nabla A(w)\bigr]
     \cdot\nabla\bigl[\mathrm{sgn}_{\eta}(A(w)-A(k))\phi\bigr]
     \nonumber\\  \qquad
     {}- \mathrm{sgn}_{\eta}(A(w)-A(k)) \mathrm{div} V(x) f(k) \phi\Bigr)\, dt\, dx =0.
  \label{weak_entropy}
\end{gather}

Note that since
$A(r)>A(k)$ if and only if $r>k$
(here we make use of the assumption that $k\in$ ``parabolic
region'', i.e., $A(k)\notin H$),
$\mathrm{sgn}_{\eta}(A(r)-A(k))\rightarrow 1$ as
$\eta\downarrow 0$ for any $r>k$. Similarly for $r<k$.
Consequently,
as $\eta\downarrow 0$,
$\int_k^w \mathrm{sgn}_{\eta}(A(\xi)-A(k))\, d\xi \to |w-k|$ a.e.~in $Q_T$.
Moreover we have
$\bigl|\int_k^w \mathrm{sgn}_{\eta}(A(\xi)-A(k))\, d\xi\bigr|
\le |w-c|\in L^1_{\mathrm{loc}}(Q_T)$ so that by
Lebesgue's dominated convergence theorem
\[
\lim_{\eta\downarrow 0} \iint\limits_{Q_T}
\biggl(\int_k^w \mathrm{sgn}_{\eta}(A(\xi)-A(k))\, d\xi
\biggr)\partial_t  \phi \, dt\, dx
= \iint\limits_{Q_T}  |w-k|\partial_t  \phi \, dt\, dx.
\]

Next we have
\begin{gather*}
   \lim_{\eta\downarrow 0}\iint\limits_{Q_T}
        \bigl[V(x)(f(w)-f(k)) - \nabla A(w)\bigr]
   \cdot\nabla \bigl[\mathrm{sgn}_{\eta}(A(w)-A(k))\phi\bigr]\, dt\, dx
   \\
\qquad {} = \lim_{\eta\downarrow 0}\iint\limits_{Q_T}
        \bigl[V(x)(f(w)-f(k)) - \nabla A(w)\bigr]
   \cdot\nabla\mathrm{sgn}_{\eta}(A(w)-A(k))\phi\, dt\, dx
   \\
\qquad {}+\lim_{\eta\downarrow 0}\iint\limits_{Q_T}
        \bigl[V(x)(f(w)-f(k)) - \nabla A(w)\bigr]
   \cdot \mathrm{sgn}_{\eta}(A(w)-A(k))\nabla\phi\, dt\, dx
   \\ \qquad {} = \underbrace{\lim_{\eta\downarrow 0}
   \iint\limits_{Q_T} V(x)(f(w)-f(k))
   \mathrm{sgn}_{\eta}'(A(w)-A(k))\nabla A(w)\phi\, dt\, dx}_{I_1}
   \\  \qquad
   {}-\lim_{\eta\downarrow 0}\iint\limits_{Q_T}
   \bigl|\nabla A(w)\bigr|^2
   \mathrm{sgn}_{\eta}'(A(w)-A(k))\phi\, dt\, dx
   \\ \qquad {}+ \underbrace{\lim_{\eta\downarrow 0}\iint\limits_{Q_T}
       \mathrm{sgn}_{\eta}(A(w)-A(k)) \bigl[V(x)(f(w)-f(k)) - \nabla A(w)\bigr]
   \cdot \nabla \phi\, dt\, dx}_{I_2}.
\end{gather*}

Note that $I_1$ can be rewritten as $I_1 = \lim\limits_{\eta\downarrow 0}\iint\limits_{Q_T}
   V(x)\mathrm{div} \mathcal{Q}_{\eta}(A(w)) \phi\, dt\, dx$, where
\begin{gather*}
   \mathcal{Q}_{\eta}(z) := \int_0^z \mathrm{sgn}_{\eta}'(r-A(k))
   \Bigl(f(A^{-1}(r)) -f(A^{-1}(A(k))) \Bigr)\,dr\\
\phantom{\mathcal{Q}_{\eta}(z) :} {}=\frac{1}{\eta}\int_{\min(z,A(k)-\eta)}^{\min(z,A(k)+\eta)}
   \Bigl(f(A^{-1}(r)) -f(A^{-1}(A(k))) \Bigr)\,dr.
\end{gather*}
Surely $\mathcal{Q}_{\eta}(z)$ tends to
zero as $\eta\downarrow 0$ for all
$z\in \text{Range(A)}$. By invoking Lebesgue's
dominated convergence theorem, we conclude after
an integration by parts that
\[
   I_1 = - \lim_{\eta\downarrow 0}\iint_{Q_T}\Bigl(
    \mathcal{Q}_{\eta}(A(w)) V(x)\cdot \nabla \phi
    + \mathcal{Q}_{\eta}(A(w)) \mathrm{div} V(x) \phi \Bigr)\, dt\, dx = 0.
\]

Using that
$\mathrm{sgn}(w-k) = \mathrm{sgn}\bigl( A(w)-A(k) \bigr)$ a.e.~in $Q_T$ (since $A(k)\notin H$)
we have
\begin{gather*}
   I_2 = \lim_{\eta\downarrow 0}\iint\limits_{Q_T}
   \mathrm{sgn}_{\eta}(A(w)-A(k)) \bigl[V(x)(f(w)-f(k)) - \nabla A(w)\bigr]
   \cdot \nabla\phi\, dt\, dx
   \\ \phantom{I_2}{}= \iint\limits_{Q_T} \mathrm{sgn}(w-k)
   \bigl[V(x)(f(w)-f(k)) - \nabla A(w)\bigr]\cdot \nabla\phi\, dt\, dx.
\end{gather*}
For the same reason we have that
\[
\lim_{\eta\downarrow 0}
\iint\limits_{Q_T} \mathrm{sgn}_{\eta}(A(w)-A(k)) \mathrm{div} V(x) f(k) \phi \, dt\, dx
=\iint\limits_{Q_T} \mathrm{sgn}(w-k) \mathrm{div} V(x) f(k) \phi \, dt\, dx.\!
\]
Consequently, letting $\eta\downarrow 0$ in \eqref{weak_entropy}, we obtain~\eqref{par_entw}.
\end{proof}

\begin{remark}
   \label{viktig_remark}
   Observe that, if $A(\cdot)$ is increasing, then a weak
   solution is automatically an entropy solution
   and hence it is unique.
\end{remark}

The next lemma, which deals with approximate
entropy inequalities for the viscous solution~$w^{\varepsilon}$, is
a direct consequence of the definition
of a weak solution of~\eqref{viscous_model}.

\begin{lemma}
   \label{viscous_hypar}
   Let $E^{\mathrm{hyp}}$ and $E^{\mathrm{par}}$ be
   defined in \eqref{Ehyp} and \eqref{Epar} respectively.
   Furthermore define
   \begin{gather}
      R_{\mathrm{visc}}:=\varepsilon \iint\limits_{Q_T}
      \bigl|\nabla w^{\varepsilon}\cdot\nabla\phi\bigr|\, dt\, dx.
   \end{gather}
   The unique weak solution $w^{\varepsilon}\in L^2\big(0,T;H^1\big({\mathbb R}^d\big)\big)\cap
   L^{\infty}\big(0,T;BV\big({\mathbb R}^d\big)\big)$ of~\eqref{viscous_model} satisfies:
   \begin{itemize}
      \item[{\rm (i)}]
      For all $k\in{\mathbb R}$ and all nonnegative $\phi\in C^{\infty}_0(Q_T)$ we have
      \begin{gather}
         \label{hyp_entw_eps}
         E^{\mathrm{hyp}}(w^{\varepsilon},k,\phi)\geq
         -R_{\mathrm{visc}}.
      \end{gather}
      We refer to \eqref{hyp_entw_eps} as
      an approximate hyperbolic entropy inequality.
      \item[{\rm (ii)}]
      For all $k\in{\mathbb R}$ such that $A(k)\notin H$ and
      all nonnegative $\phi\in C^{\infty}_0(Q_T)$ we have
      \begin{gather}
         \label{par_entw_eps}
         E^{\mathrm{par}}(w^{\varepsilon},k,\phi)\geq -R_{\mathrm{visc}}.
      \end{gather}
      We refer to \eqref{par_entw_eps} as
     an approximate parabolic entropy inequality.
   \end{itemize}
\end{lemma}

\begin{proof}
In what follows we always let $k$ and $\phi$ be as indicated by the lemma.
The proof of the inequality~\eqref{hyp_entw_eps} follows
the proof of~\eqref{hyp_entw} rather closely.
Since $w^{\varepsilon}$ is a weak solution and
$\bigl[\mathrm{sgn}_{\eta}(w^{\varepsilon}-k)\phi\bigr]$ belongs to
$L^2\big(0,T; H^1\big({\mathbb R}^d\big)\big)$, we have
\begin{gather*}
    -\int_0^T \Bigl\langle \partial_t  w^{\varepsilon},
    \mathrm{sgn}_{\eta}(w^{\varepsilon}-k)\phi \Bigr\rangle \, dt\\
\qquad {}    + \iint\limits_{Q_T} \Bigl( \bigl[V(x)(f(w^{\varepsilon})-f(k))
    - \nabla A^{\varepsilon}(w^{\varepsilon}) \bigr]
    \cdot\nabla\bigl[\mathrm{sgn}_{\eta}(w^{\varepsilon}-k)\phi\bigr]
   \\  \qquad
   {}- \mathrm{div} V(x) f(k)\bigl[\mathrm{sgn}_{\eta}(w^{\varepsilon}-k)\phi\bigr]\Bigr)\, dt\, dx =0.
\end{gather*}
By the chain rule we obviously have
\begin{gather*}
   -\int_0^T \Bigl\langle \partial_t  w^{\varepsilon},
    \mathrm{sgn}_{\eta}(w^{\varepsilon}-k)\phi\Bigr\rangle \, dt
    \\  \qquad {} = \iint_{Q_T}
    \biggl(\int_k^{w^{\varepsilon}} \mathrm{sgn}_{\eta}(\xi-k)\, d\xi \biggr)
    \partial_t \phi\, dt\, dx
    \overset{\eta\downarrow 0}{\longrightarrow}
    \iint\limits_{Q_T}  |w^{\varepsilon}-k|\partial_t  \phi \, dt\, dx
\end{gather*}
so that
\begin{gather}
     \iint\limits_{Q_T} |w^{\varepsilon}-k|\partial_t  \phi \, dt\, dx \nonumber\\
\qquad {}     +\lim_{\eta\downarrow 0}\iint\limits_{Q_T}
     \Bigl(\bigl[V(x)(f(w^{\varepsilon})-f(k))
     - \nabla A^{\varepsilon}(w^{\varepsilon})\bigr]
     \cdot\nabla\bigl[\mathrm{sgn}_{\eta}(w^{\varepsilon}-k)\phi\bigr]
     \nonumber\\  \qquad
    {} - \mathrm{sgn}_{\eta}(w^{\varepsilon}-k) \mathrm{div} V(x) f(k) \phi\Bigr)\, dt\, dx =0.
  \label{weak_entropy_eps}
\end{gather}
Firstly we have
\[
\lim_{\eta\downarrow 0}
\iint\limits_{Q_T} \mathrm{sgn}_{\eta}(w^{\varepsilon}-k) \mathrm{div} V(x) f(k) \phi \, dt\, dx
=\iint\limits_{Q_T} \mathrm{sgn}(w^{\varepsilon}-k) \mathrm{div} V(x) f(k) \phi \, dt\, dx.
\]
Next we have
\begin{gather*}
   \lim_{\eta\downarrow 0}\iint\limits_{Q_T}
        \bigl[V(x)(f(w^{\varepsilon})-f(k)) - \nabla A^{\varepsilon}(w^{\varepsilon})\bigr]
   \cdot\nabla \bigl[\mathrm{sgn}_{\eta}(w^{\varepsilon}-k)\phi\bigr]\, dt\, dx
   \\ \qquad {}= \lim_{\eta\downarrow 0}\iint\limits_{Q_T}
        \bigl[V(x)(f(w^{\varepsilon})-f(k)) - \nabla A^{\varepsilon}(w^{\varepsilon})\bigr]
   \cdot\nabla\mathrm{sgn}_{\eta}(w^{\varepsilon}-k)\phi\, dt\, dx
   \\ \qquad {}+\lim_{\eta\downarrow 0}\iint\limits_{Q_T}
        \bigl[V(x)(f(w^{\varepsilon})-f(k)) - \nabla A^{\varepsilon}(w^{\varepsilon})\bigr]
   \cdot \mathrm{sgn}_{\eta}(w^{\varepsilon}-k)\nabla\phi\, dt\, dx
   \\ \qquad{} = \underbrace{\lim_{\eta\downarrow 0}
   \iint\limits_{Q_T} V(x)(f(w^{\varepsilon})-f(k))
   \mathrm{sgn}_{\eta}'
(w^{\varepsilon}-k)\nabla A^{\varepsilon}(w^{\varepsilon})\phi\, dt\, dx}_{I_1}
   \\ \qquad
   {}-\lim_{\eta\downarrow 0}\iint\limits_{Q_T} (A^{\varepsilon})'(w^{\varepsilon})
   \bigl|\nabla w^{\varepsilon}\bigr|^2
   \mathrm{sgn}_{\eta}'(w^{\varepsilon}-k)\phi\, dt\, dx
   \\  \qquad{} + \iint\limits_{Q_T} \mathrm{sgn}(w^{\varepsilon}-k)
       \bigl[V(x)(f(w^{\varepsilon})-f(k)) - \nabla A^{\varepsilon}(w^{\varepsilon})\bigr]
   \cdot \nabla \phi\, dt\, dx.
\end{gather*}

Note that $I_1$ can be rewritten as $I_1 = \lim\limits_{\eta\downarrow 0}\iint\limits_{Q_T}
   V(x)\mathrm{div} \mathcal{Q}_{\eta}(w^{\varepsilon}) \phi\, dt\, dx$, where
\begin{gather*}
   \mathcal{Q}_{\eta}(z) := \int_0^z \mathrm{sgn}_{\eta}'(r-k)
   \bigl(f(r) -f(k) \bigr)\,dr\\
\phantom{\mathcal{Q}_{\eta}(z) :} {}   =\frac{1}{\eta}\int_{\min(z,k-\eta)}^{\min(z,k+\eta)}
   \bigl(f(r) -f(k) \bigr)\,dr \to 0 \,\, \text{as $\eta\downarrow 0$}.
\end{gather*}

 From Lebesgue's
dominated convergence theorem we conclude that
\[
   I_1 = - \lim_{\eta\downarrow 0}\iint_{Q_T}\Bigl(
    \mathcal{Q}_{\eta}(A(w^{\varepsilon})) V(x)\cdot \nabla \phi
    + \mathcal{Q}_{\eta}(A(w^{\varepsilon})) \mathrm{div} V(x) \phi \Bigr)\, dt\, dx = 0.
\]

In conclusion we have
\begin{gather}
\iint\limits_{Q_T}\Bigl( |w^{\varepsilon}-k|\phi_t
      + \mathrm{sgn}(w^{\varepsilon}-k)\bigl[V(x)(f(w^{\varepsilon})- f(k))
      - \nabla A^{\varepsilon}(w^{\varepsilon})\bigr]\cdot\nabla\phi
      \nonumber\\ \qquad
    {}  -\mathrm{sgn}(w^{\varepsilon}-k)\mathrm{div} V(x) f(k)\phi\Bigr)\, dt\, dx
      \nonumber\\  \qquad
      {}= \lim_{\varepsilon\downarrow 0} \iint\limits_{Q_T}
      (A^{\varepsilon})'(w^{\varepsilon})\bigl|\nabla w^{\varepsilon}\bigr|^2
      \mathrm{sgn}_{\eta}(w^{\varepsilon}-k)\phi\, dt\, dx\ge 0
   \label{eqn:Mod_Carrillo}
\end{gather}
for any $0\le\phi\in C^{\infty}_0(Q_T)$ and any $k\in{\mathbb R}$.
 From this we conclude easily that \eqref{hyp_entw_eps} holds.

It remains to prove the parabolic entropy inequality \eqref{par_entw_eps}.
Let $0\le\phi\in C^{\infty}_0(Q_T)$ and $k\in{\mathbb R}$ be
such that $A(k)\notin H$.
Starting off by choosing $\bigl[\mathrm{sgn}_{\eta}(A(w^{\varepsilon})-A(k))\phi\bigr]$
as a test function in the weak formulation and
then continuing exactly as in the proof of \eqref{par_entw}, we obtain
\begin{gather*}
   E^{\mathrm{par}}(w^{\varepsilon},k,\phi) =
   \lim_{\eta\downarrow 0}\iint\limits_{Q_T}
   \varepsilon \nabla w^{\varepsilon}\cdot \nabla \bigl[\mathrm{sgn}_{\eta}
   \bigl(A(w^{\varepsilon})-A(k)\bigr)\phi\bigr]\, dt\, dx.
\end{gather*}
The right-hand side of this equality can be expanded into
\begin{gather*}
\lim_{\eta\downarrow 0} \iint\limits_{Q_T}\Bigl(
   \varepsilon A'(w^{\varepsilon})\bigl|\nabla w^{\varepsilon}\bigr|^2
   \mathrm{sgn}_{\eta}'(A(w^{\varepsilon})-A(k))\phi \\
\qquad {}   +  \varepsilon\, \mathrm{sgn}_{\eta}(A(w^{\varepsilon})-A(k))\nabla w^{\varepsilon}
   \cdot\nabla\phi\Bigr)\, dt\, dx \nonumber\\
 \qquad {}\geq  \lim_{\eta\downarrow 0}\iint\limits_{Q_T} \varepsilon\,
   \mathrm{sgn}_{\eta}(A(w^{\varepsilon})-A(k)) \nabla w^{\varepsilon} \cdot\nabla\phi\, dt\, dx
   \ge  - \varepsilon \iint\limits_{Q_T} \bigl
  |\nabla w^{\varepsilon}\cdot \nabla\phi\bigr|\, dt\, dx.
\end{gather*}
This concludes the proof of \eqref{par_entw_eps}.
\end{proof}

\section{Proof of Theorem \ref{main_thm}}
\label{sec:main_thm}
Following Carrillo \cite{Carrillo} (see also \cite{KR})
in this section we use Lemma \ref{entropy_hypar} and
Lemma~\ref{viscous_hypar} to prove Theorem~\ref{main_thm}.
Let $w^{\varepsilon}=w^{\varepsilon}(x,t)$ solve~\eqref{cont_model}
and $w=w(y,s)$ solve~\eqref{viscous_model}.
Following Kru\v{z}kov~\cite{Kruzkov} and
Kuznetsov~\cite{Kuznetsov:1976} we now specify a nonnegative test
function $\phi=\phi(t,x,s,y)$ defined on $Q_T\times Q_T$. To this
end let $\rho\in C_0^{\infty}({\mathbb R})$ be a function satisfying
\[
\mathrm{supp}(\rho)\subset \{\sigma\in{\mathbb R}:|\sigma|\leq 1\},\qquad
\rho(\sigma)\geq 0\,\forall \sigma\in{\mathbb R}, \qquad
\int_{{\mathbb R}}\rho(\sigma)\, d\sigma=1.
\]
For $x\in{\mathbb R}^d,t\in{\mathbb R}$ and $r,r_0>0$, let $\omega_{r}(x)=
\frac{1}{r}\rho\left(\frac{x_1}{r}\right)\cdots
\frac{1}{r}\rho\left(\frac{x_d}{r}\right)$ and
$\rho_{r_0}(t)=\frac{1}{r_0}\rho\left(\frac{t}{r_0}\right)$.
Pick any two points $\nu,\tau\in (0,T)$, $\nu<\tau$.
For any $\alpha_0>0$ define
\[
\psi_{\alpha_0}(t)=H_{\alpha_0}(t-\nu) - H_{\alpha_0}(t-\tau),
\qquad H_{\alpha_0}(t)=\int_{-\infty}^{t}
\rho_{\alpha_0}(\xi)\, d\xi.
\]
With $0<r_0 < \min(\nu,T-\tau)$ and
$\alpha_0\in \bigl(0,\min(\nu-r_0,T-\tau-r_0)\bigr)$ we set
\begin{equation}
   \label{test}
   \phi(x,t,y,s)
   := \psi_{\alpha_0}(t)\omega_r(x-y)\rho_{r_0}(t-s).
\end{equation}
Note that $\mathrm{supp}(\phi(x,\cdot,y,s))\subset (r_0,T-r_0)$
for all $x,y\in {\mathbb R}^d,s\in (0,T)$ and
$\mathrm{supp}(\phi(x,t,y,\cdot)) $ $\subset (0,T)$
for all $x,y\in {\mathbb R}^d,t\in (0,T)$. Consequently
$(x,t)\mapsto \phi(x,t,y,s)$ belongs to $C^{\infty}_0(Q_T)$ for
each fixed $(y,s)\in Q_T$ and $(y,s)\mapsto \phi(x,t,y,s)$ belongs to
$C^{\infty}_0(Q_T)$ for each fixed $(x,t)\in Q_T$.

Observe that with the choice of $\phi$
as in \eqref{test} we have
\begin{gather}
      \partial_t \phi + \partial_s \phi=
   \bigl[\rho_{\alpha_0}(t-\nu)- \rho_{\alpha_0}(t-\tau)\bigr]
   \omega_r(x-y)\rho_{r_0}(t-s),\nonumber\\
   \nabla_x\phi + \nabla_y\phi=0.\label{phi_prop}
\end{gather}
Before continuing we need to introduce the two ``hyperbolic'' sets
\begin{gather*}
   \mathcal{H}^{\varepsilon} = \Bigl\{(x,t)\in Q_T: A(w^{\varepsilon}(x,t))\in H \Bigr\},
   \qquad
   \mathcal{H} = \Bigl\{(y,s)\in Q_T: A(w(y,s))\in H \Bigr\}
\end{gather*}
and note that
\begin{gather}
   \label{hyp_set}
   \nabla_x A(w^{\varepsilon})=0\;\; \text{a.e.~in $\mathcal{H}^{\varepsilon}$}
   \quad \text{and} \quad
   \nabla_y A(w)=0\;\;\text{a.e.~in $\mathcal{H}$}, \\
   \mathrm{sgn}(w^{\varepsilon}-w)  =\mathrm{sgn}\bigl(A(w^{\varepsilon})-A(w)\bigr) \nonumber\\
\qquad \text{a.e.~in
   $\Bigl[(Q_T\setminus\mathcal{H})\times Q_T \Bigr]
   \cup \Bigl[Q_T\times (Q_T\setminus\mathcal{H}^{\varepsilon})\Bigr]$}.
   \label{par_set}
\end{gather}

Using the approximate hyperbolic entropy
inequality~\eqref{hyp_entw_eps}
for the viscous solution $w^{\varepsilon}=w^{\varepsilon}(x,t)$ with $k=w(y,s)$, we
get for $(y,s)\in Q_T$
\begin{gather}
\iint\limits_{Q_T}
      \Bigl(  |w^{\varepsilon}-w|\partial_t  \phi
      + \mathrm{sgn}(w^{\varepsilon}-w)\bigl[V(x)(f(w^{\varepsilon})- f(w))
      - \nabla_x A(w^{\varepsilon})\bigr]
      \cdot\nabla_x\phi
     \nonumber\\  \qquad
    {} - \mathrm{sgn}(w^{\varepsilon}-w)\mathrm{div}_x V(x)f(w)\phi\Bigr)\, dt\, dx\, ds\, dy
     \geq - \overline{R}_{\mathrm{visc}}.   \label{v_entropy_tmpI}
\end{gather}
Using the approximate parabolic entropy
inequality~\eqref{par_entw_eps}
for the viscous solution $w^{\varepsilon}=w^{\varepsilon}(x,t)$ with $k=w(y,s)$, we
get for $(y,s)\in Q_T\setminus\mathcal{H}$
\begin{gather}
\iint\limits_{Q_T}
      \Bigl(  |w^{\varepsilon}-w|\partial_t  \phi
      + \mathrm{sgn}(w^{\varepsilon}-w)\bigl[V(x)(f(w^{\varepsilon})- f(w))
      - \nabla_x A(w^{\varepsilon})\bigr]
      \cdot\nabla_x\phi
    \nonumber \\  \qquad
{}     - \mathrm{sgn}(w^{\varepsilon}-w)\mathrm{div}_x V(x)f(w)\phi\Bigr)\, dt\, dx
     \nonumber\\  \qquad {}\geq
     \lim_{\eta\downarrow 0}\iint\limits_{Q_T}
     \bigl|\nabla_x A(w^{\varepsilon})\bigr|^2
     \mathrm{sgn}_{\eta}'(A(w^{\varepsilon})-A(w))\phi\, dt\, dx
     - \overline{R}_{\mathrm{visc}}.   \label{v_entropy_tmpII}
\end{gather}

Next we would like to integrate \eqref{v_entropy_tmpI}
and \eqref{v_entropy_tmpII} over $(y,s)\in Q_T$ and
$(y,s)\in Q_T\setminus\mathcal{H}$ respectively.
To this end we need to know that the involved functions
are $(y,s)$ integrable. Consider first
$(y,s)\mapsto \iint\limits_{Q_T}\mathrm{sgn}(v-u)\nabla_x A(w^{\varepsilon})
\cdot \nabla_x\phi\, dt\, dx$.
We denote this function by~$D(y,s)$.

To see that $D(\cdot,\cdot)$ is integrable on $Q_T$ we observe that
for each fixed $(y,s)\in Q_T$
\[
\mathrm{sgn}(v-u)\nabla_x A(w^{\varepsilon}) = \nabla_x \left|A(w^{\varepsilon})-A(w)\right|
\;\; \text{for a.e.~$(x,t)\in Q_T$}
\]
and hence
\begin{gather*}
   D(y,s) =\iint\limits_{Q_T}
   \Bigl[\nabla_x\bigl|A(w^{\varepsilon})-A(w)\bigr|\Bigr]
   \cdot \nabla_x \phi \, dt\, dx.
\end{gather*}
Since the function $(x,t)\mapsto \phi(x,t,y,s)$
belongs to $C^{\infty}_0(Q_T)$ for each fixed
$(y,s)\in Q_T$, an integration by parts in $x$ gives
\begin{gather*}
   D(y,s) =- \iint\limits_{Q_T}
   \bigl|A(w^{\varepsilon})-A(w)\bigr| \Delta_x \phi\, dt\, dx.
\end{gather*}
Integration over $(y,s)\in Q_T$ and estimation yield
\begin{gather*}
   \biggl|\;\iint\limits_{Q_T}D(y,s)\, ds\, dy\biggr|
   \le \iiiint\limits_{Q_T\times Q_T}
   \Bigl( |A(w^{\varepsilon}(x,t))| + |A(w(y,s))|\Bigr)
   \Delta_x \phi(x,y,t,s)\, dt\, dx\, ds\, dy.
\end{gather*}
By changing the variables ($z:=x-y$, $\tau=t-s$) and
taking into account that $w^{\varepsilon},w\in L^1(Q_T)$ we find that
\begin{gather*}
   \biggl|\;\iint\limits_{Q_T} D(y,s)\, ds\, dy\biggr|
   \le \iiiint
   |A(w^{\varepsilon}(x,t))|\psi_{\alpha_0}(t)
   |\Delta_z \omega_r(z)|\, \rho_{r_0}(\tau)\, dt\, dx\, d\tau\, dz
    \\  \qquad {} + \iiiint |A(w(x-z,t-\tau))| \psi_{\alpha_0}(t)
   |\Delta_z \omega_r(z)|\, \rho_{r_0}(\tau)\, dt\, dx\, d\tau\, dz
   \\ \qquad{} \le \|A(w^{\varepsilon})\|_{L^1(Q_T)}
    \| \Delta_z\omega_r \|_{L^1({\mathbb R}^d)}
   + \|A(w)\|_{L^1(Q_T)}\|\Delta_z \omega_r \|_{L^1({\mathbb R}^d)} <\infty.
\end{gather*}
Hence we have that $D(\cdot,\cdot)$ is integrable on $Q_T$.

In a similar vein one can also show the
integrability of
\begin{gather*}
(y,s) \mapsto \iint\limits_{Q_T} |w^{\varepsilon}-w|\partial_t\phi\, dt\, dx, \\
(y,s) \mapsto \iint\limits_{Q_T} \mathrm{sgn}(w^{\varepsilon}-w)
V(x)(f(w^{\varepsilon})- f(w))\cdot \nabla_x \phi \, dt\, dx, \\
(y,s) \mapsto \iint\limits_{Q_T} \mathrm{sgn}(w^{\varepsilon}-w)
\mathrm{div}_x V(x)f(w)\phi\, dt\, dx, \qquad {\rm and} \quad
(y,s) \mapsto \overline{R}_{\mathrm{visc}}.
\end{gather*}

It remains to consider the
integrability of the function
\[
Q_T\setminus\mathcal{H}\ni (y,s) \mapsto \lim_{\eta\downarrow 0}\iint\limits_{Q_T}
     \bigl|\nabla_x A(w^{\varepsilon})\bigr|^2
     \mathrm{sgn}_{\eta}'(A(w^{\varepsilon})-A(w))\phi\, dt\, dx.\]
  This follows
from \eqref{v_entropy_tmpII}. We have
by Lebesgue's dominated convergence theorem and
the first part of~\eqref{hyp_set}
\begin{gather}
\iint\limits_{Q_T\setminus\mathcal{H}}
      \biggl(\lim_{\eta\downarrow 0}
      \iint\limits_{Q_T}
      \bigl|\nabla_x A(w^{\varepsilon})\bigr|^2
      \mathrm{sgn}_{\eta}'(A(w^{\varepsilon})-A(w))\phi\, dt\, dx\biggr)\, ds\, dy
      \nonumber\\ \qquad
     {} = \lim_{\eta\downarrow 0}
      \iiiint\limits_{(Q_T\setminus\mathcal{H})\times Q_T}
      \bigl|\nabla_x A(w^{\varepsilon})\bigr|^2
      \mathrm{sgn}_{\eta}'(A(w^{\varepsilon})-A(w))\phi\, dt\, dx\, ds\, dy.
      \nonumber\\  \qquad {}= \lim_{\eta\downarrow 0}
      \iiiint\limits_{(Q_T\setminus\mathcal{H})\times (Q_T\setminus\mathcal{H})}
      \bigl|\nabla_x A(w^{\varepsilon})\bigr|^2
      \mathrm{sgn}_{\eta}'(A(w^{\varepsilon})-A(w))\phi\, dt\, dx\, ds\, dy.
   \label{Diss}
\end{gather}

We now integrate \eqref{v_entropy_tmpI} over $(y,s)\in Q_T$
and \eqref{v_entropy_tmpII} over $(y,s)\in Q_T\setminus\mathcal{H}$.
Addition of the two resulting inequalities yields
\begin{gather}
\iiiint\limits_{Q_T\times Q_T}
      \Bigl(  |w^{\varepsilon}-w|\partial_t  \phi
      + \mathrm{sgn}(w^{\varepsilon}-w)\bigl[V(x)(f(w^{\varepsilon})- f(w))
      - \nabla_x A(w^{\varepsilon})\bigr]
      \cdot\nabla_x\phi
     \nonumber\\  \qquad
     {}- \mathrm{sgn}(w^{\varepsilon}-w)\mathrm{div}_x V(x)f(w)\phi\Bigr)\, dt\, dx\, ds\, dy
     \nonumber\\  \qquad{}= \!\!\iiiint\limits_{(Q_T\setminus\mathcal{H})\times Q_T}\!\!
     \Bigl(|w^{\varepsilon}-w|\partial_t  \phi
     +\mathrm{sgn}(w^{\varepsilon}-w)\bigl[V(x)(f(w^{\varepsilon})- f(w))
     - \nabla_x A(w^{\varepsilon})\bigr]\cdot\nabla_x\phi
    \nonumber \\  \qquad
{}- \mathrm{sgn}(w^{\varepsilon}-w)\mathrm{div}_xV(x)f(w)\phi\Bigr)\, dt\, dx\, ds\, dy
     \nonumber\\  \qquad {} + \iiiint\limits_{\mathcal{H}\times Q_T}
     \Bigl (|w^{\varepsilon}-w|\partial_t  \phi
     +\mathrm{sgn}(w^{\varepsilon}-w)
\bigl[V(x)(f(w^{\varepsilon})- f(w))- \nabla_x A(w^{\varepsilon})\bigr]
     \cdot\nabla_x\phi
    \nonumber \\  \qquad
    {} - \mathrm{sgn}(w^{\varepsilon}-w)\mathrm{div}_x V(x)f(w)\phi\Bigr)\, dt\, dx\, ds\, dy
     \nonumber\\  \qquad {}\ge  \lim_{\eta\downarrow 0}
     \iiiint\limits_{(Q_T\setminus\mathcal{H})
     \times (Q_T\setminus\mathcal{H}^{\varepsilon})}
     \bigl|\nabla_x A(w^{\varepsilon})\bigr|^2
     \mathrm{sgn}_{\eta}'(A(w^{\varepsilon})-A(w))\phi\, dt\, dx\, ds\, dy
     - \overline{R}_{\mathrm{visc}},\!\!\!
   \label{v_entropy}
\end{gather}
where $\overline{R}_{\mathrm{visc}}:=
\iint\limits_{Q_T} R_{\mathrm{visc}}\, ds\, dy$
and we have used~\eqref{Diss}.

Similarly, using the hyperbolic, parabolic entropy
inequalities \eqref{hyp_entw}, \eqref{par_entw} for
the exact entropy solution $w=w(y,s)$ with $k=w^{\varepsilon}(x,t)$
and then integrating over $(x,t)\in Q_T$, we get
\begin{gather}
  \iiiint\limits_{Q_T\times Q_T}
      \Bigl( |w-w^{\varepsilon}|\partial_s \phi
      + \mathrm{sgn}(w-w^{\varepsilon})\bigl[V(y)(f(w)- f(w^{\varepsilon}))
      - \nabla_y A(w)\bigr]\cdot\nabla_y\phi
      \nonumber\\  \qquad
      {}- \mathrm{sgn}(w-w^{\varepsilon})\mathrm{div}_y V(y)f(w^{\varepsilon})\phi\Bigr)\, dt\, dx\, ds\, dy
      \nonumber\\
\qquad {}\geq \lim_{\eta\downarrow 0}
      \iiiint\limits_{(Q_T\setminus\mathcal{H}^{\varepsilon})
      \times(Q_T\setminus\mathcal{H})}
      \bigl|\nabla_yA(w)\bigr|^2 \mathrm{sgn}_{\eta}'(A(w)-A(w^{\varepsilon}))\phi\, dt\, dx\, ds\, dy.
 \label{w_entropy}
\end{gather}

Using \eqref{hyp_set} and \eqref{par_set} we find that
\begin{gather}
 \iiiint\limits_{Q_T\times Q_T} \mathrm{sgn}(w_{\varepsilon}-w)
        \nabla_x A(w^{\varepsilon}) \cdot \nabla_y\phi\, dt\, dx\, ds\, dy
     \nonumber \\\qquad
     {} = - \iiiint\limits_{Q_T\times (Q_T\setminus\mathcal{H}^{\varepsilon})}
        \mathrm{sgn}(A(w^{\varepsilon}) - A(w))\nabla_x A(w^{\varepsilon})
      \cdot\nabla_y\phi\, dt\, dx\, ds\, dy
      \nonumber\\ \qquad {}= - \lim_{\eta\downarrow 0}\iiiint
      \limits_{Q_T\times (Q_T\setminus\mathcal{H}^{\varepsilon})}
      \mathrm{sgn}_{\eta}(A(w^{\varepsilon}) - A(w))\nabla_x A(w^{\varepsilon})
      \cdot\nabla_y\phi\, dt\, dx\, ds\, dy
      \nonumber\\ \qquad {}= - \lim_{\eta\downarrow 0}\iiiint
      \limits_{Q_T\times (Q_T\setminus\mathcal{H}^{\varepsilon})}
       \nabla_yA(w)\cdot\nabla_x A(w^{\varepsilon})
      \mathrm{sgn}_{\eta}'(A(w^{\varepsilon}) - A(w))\phi\, dt\, dx\, ds\, dy.
      \nonumber\\ \qquad {}= - \lim_{\eta\downarrow 0}\!\!\iiiint
      \limits_{(Q_T\setminus\mathcal{H})
         \times (Q_T\setminus\mathcal{H}^{\varepsilon}) } \!\!\!\!\!
       \nabla_yA(w)\cdot\nabla_x A(w^{\varepsilon})
       \mathrm{sgn}_{\eta}'(A(w^{\varepsilon})\! - A(w))\phi\, dt\, dx\, ds\, dy.
   \label{cross_v_ny}
\end{gather}

Similarly, again
using \eqref{hyp_set} and \eqref{par_set}, we find that
\begin{gather}
- \iiiint\limits_{Q_T\times Q_T} \mathrm{sgn} (w-w^{\varepsilon})\nabla_y A(w)
      \cdot\nabla_x\phi\, dt\, dx\, ds\, dy
      \nonumber\\ \qquad{}= - \lim_{\eta\downarrow 0}\!\!\iiiint
      \limits_{(Q_T\setminus\mathcal{H}) \times
      (Q_T\setminus\mathcal{H}^{\varepsilon})} \!\!\!\!\!
      \nabla_x A(w^{\varepsilon})\cdot \nabla_y A(w)
      \mathrm{sgn}_{\eta}'(A(w)\!-A(w^{\varepsilon}))\phi\, dt\, dx\, ds\, dy.\!
   \label{cross_w_ny}
\end{gather}

The use of the second part of \eqref{phi_prop} when adding
\eqref{v_entropy} and \eqref{cross_v_ny} yields
\begin{gather}
\iiiint\limits_{Q_T\times Q_T} \Bigl(  |w^{\varepsilon}-w|\partial_t  \phi
      + \mathrm{sgn}(w^{\varepsilon}-w) \bigl[ V(x)(f(w^{\varepsilon})- f(w))\bigr]\cdot\nabla_x \phi
     \nonumber \\  \qquad {}
      - \mathrm{sgn}(w^{\varepsilon}-w)\mathrm{div}_xV(x)f(w)\phi \Bigr) \, dt\, dx\, ds\, dy
      \nonumber\\ \qquad
      {}\ge \lim_{\eta\downarrow 0}
      \iiiint\limits_{(Q_T\setminus\mathcal{H}^{\varepsilon})
      \times (Q_T\setminus\mathcal{H})}
      \Bigl( \bigl|\nabla_x A(w^{\varepsilon})\bigr|^2 -
      \nabla_yA(w) \cdot\nabla_x A(w^{\varepsilon})\Bigr)
      \nonumber\\ \qquad
      {}\times \mathrm{sgn}_{\eta}'(A(w^{\varepsilon})-A(w))
      \phi\, dt\, dx\, ds\, dy - R_{\mathrm{visc}}.
   \label{v_entropy_ny}
\end{gather}

Similarly the addition of \eqref{w_entropy} and \eqref{cross_w_ny} yields
\begin{gather}
\iiiint\limits_{Q_T\times Q_T}\Bigl(|w-w^{\varepsilon}|\partial_s \phi
      + \mathrm{sgn}(w-w^{\varepsilon})\bigl[V(y)(f(w)- f(w^{\varepsilon})\bigr]\cdot\nabla_y \phi
     \nonumber \\  \qquad {}
      - \mathrm{sgn}(w-w^{\varepsilon})\mathrm{div}_y V(y)f(w^{\varepsilon})\phi \Bigr)\, dt\, dx\, ds\, dy
      \nonumber\\ \qquad
      {}\ge  \lim_{\eta\downarrow 0}
      \iiiint\limits_{(Q_T\setminus\mathcal{H}^{\varepsilon})
      \times (Q_T\setminus\mathcal{H})}
      \Bigl( \bigl|\nabla_yA(w)\bigr|^2 -
      \nabla_xA(w^{\varepsilon})\cdot\nabla_yA(w) \Bigr)
      \nonumber\\  \qquad {}      \times \mathrm{sgn}_{\eta}'(A(w)-A(w^{\varepsilon}))\phi\, dt\, dx\, ds\, dy.
   \label{w_entropy_ny}
\end{gather}

Following Karlsen and Risebro \cite{KR} we write
\begin{gather*}
\mathrm{sgn}(w^{\varepsilon}-w)V(x)(f(w^{\varepsilon}) - f(w))\cdot\nabla_x \phi
   -\mathrm{sgn}(w^{\varepsilon}-w) \mathrm{div}_x V(x)f(w)\phi
   \\  \qquad
   {}= \mathrm{sgn}(w^{\varepsilon}-w)\bigl(V(x)f(w^{\varepsilon}) - V(y)f(w)\bigr)\cdot\nabla_x \phi
   \\
\qquad{}
+ \mathrm{sgn}(w^{\varepsilon}-w)\mathrm{div}_x\bigl[\bigl(V(y)f(w) - V(x)f(w)\bigr)\phi\bigr],
\\ \mathrm{sgn}(w-w^{\varepsilon})V(y)(f(w)-f(w^{\varepsilon}))\cdot \nabla_y \phi
    - \mathrm{sgn}(w-w^{\varepsilon}) \mathrm{div}_y V(y)f(w^{\varepsilon})\phi
   \\  \qquad
   {}= \mathrm{sgn}(w^{\varepsilon}-w)\bigl(V(x)f(w^{\varepsilon}) - V(y)f(w)\bigr)\cdot \nabla_y \phi
   \\
\qquad {}
- \mathrm{sgn}(w^{\varepsilon}-w)\mathrm{div}_y\bigl[\bigl(V(x)f(w^{\varepsilon})
   -V(y)f(w^{\varepsilon})\bigr)\phi\bigr].
\end{gather*}

When adding \eqref{v_entropy_ny}
and \eqref{w_entropy_ny}, we use the second part
of \eqref{phi_prop} and the identities
\[
\mathrm{sgn}(-r)=-\mathrm{sgn}(r)\,\, \text{a.e.~in ${\mathbb R}$},\qquad
\mathrm{sgn}_{\eta}'(-r)=\mathrm{sgn}_{\eta}'(r)\,\, \text{a.e.~in ${\mathbb R}$}.
\]
The final result takes the form
\begin{gather}
- \iiiint\limits_{Q_T\times Q_T}
      |w^{\varepsilon}-w| (\partial_t \phi+\partial_s  \phi)\, dt\, dx\, ds\, dy
     \nonumber\\
\qquad {} \le R_{\mathrm{diss}} +  \overline{R}_{\mathrm{visc}}
      + R_{\mathrm{conv}} \le \overline{R}_{\mathrm{visc}}
      +  R_{\mathrm{conv}},   \label{total_entropy}
\end{gather}
where the expression for $\partial_t \phi+\partial_s \phi$
is written out in~\eqref{phi_prop},
\begin{gather*}
R_{\mathrm{conv}}:=\iiiint\limits_{Q_T\times Q_T}
I_{\mathrm{conv}}\, dt\, dx\, ds\, dy,\\
   I_{\mathrm{conv}}  := \mathrm{sgn}(w^{\varepsilon}-w)\Bigl(\mathrm{div}_x \bigl[\bigl(V(y)f(w)
   -V(x)f(w)\bigr)\phi\bigr]\\
\phantom{I_{\mathrm{conv}}  :=}{}
   - \mathrm{div}_y \bigl[\bigl(V(x)f(w^{\varepsilon})
   - V(y)f(w^{\varepsilon})\bigr)\phi\bigr]\Bigr),
\end{gather*}
and
\begin{gather*}
R_{\mathrm{diss}} := -\lim\limits_{\eta\downarrow  0}
      \iiiint\limits_{(Q_T\setminus\mathcal{H}^{\varepsilon})
      \times (Q_T\setminus\mathcal{H})}
      \bigl|\nabla_xA(w^{\varepsilon}) - \nabla_yA(w) \bigr|^2  \\
\phantom{R_{\mathrm{diss}} :=}{}\times      \mathrm{sgn}_{\eta}' \bigl(A(w^{\varepsilon}) - A(w)\bigr)\phi\, dt\, dx\, ds\, dy\le 0.
\end{gather*}

Having in mind the first part of
 \eqref{phi_prop}, we get by the triangle inequality
\[
-\iiiint\limits_{Q_T\times Q_T}|w^{\varepsilon}(x,t)-w(y,s)|
(\partial_t \phi+\partial_s  \phi)\, dt\, dx\, ds\, dy
\le R_{w^{\varepsilon},w} + R_{w,x} + R_{w,t},
\]
where
\begin{gather*}
   R_{w^{\varepsilon},w}:=-\iiiint\limits_{Q_T\times Q_T}|w^{\varepsilon}(x,t)-w(x,t)|
   \bigl[\rho_{\alpha_0}(t-\nu)-
   \rho_{\alpha_0}(t-\tau)\bigr]\\
\phantom{R_{w^{\varepsilon},w}:=}{} \times\omega_r(x-y)\rho_{r_0}(t-s)
   \, dt\, dx\, ds\, dy,\\
      R_{w,x}:=- \iiiint\limits_{Q_T\times Q_T}|w(x,t)-w(y,t)|
   \bigl[\rho_{\alpha_0}(t-\nu)-\rho_{\alpha_0}(t-\tau)\bigr]\\
\phantom{R_{w,x}:=}{}\times
   \omega_r(x-y)\rho_{r_0}(t-s)\, dt\, dx\, ds\, dy,
   \end{gather*}
\begin{gather*}
   R_{w,t}: = -\iiiint\limits_{Q_T\times Q_T}|w(y,t)-w(y,s)|
   \bigl[\rho_{\alpha_0}(t-\nu)-\rho_{\alpha_0}(t-\tau)\bigr]\\
\phantom{R_{w,t}: =}{}\times
   \omega_r(x-y)\rho_{r_0}(t-s)\, dt\, dx\, ds\, dy.
\end{gather*}

Firstly a standard
$L^1$ continuity argument gives
$\lim\limits_{r_0\downarrow 0} R_{w,t}=0$. Next
\begin{gather*}
   \lim_{\alpha_0\downarrow 0} R_{w,x}=
   \int_{{\mathbb R}^d}\int_{{\mathbb R}^d}\Bigl( |w(x,\tau)-w(y,\tau)| -
   |w(x,\nu)-w(y,\nu)|\Bigr)\omega_r(x-y)\, dx\, dy
   \\
\phantom{\lim_{\alpha_0\downarrow 0} R_{w,x}=} {} \overset{(z:=x-y)}{\le}
   \int_{{\mathbb R}^d}\int_{{\mathbb R}^d} |w(y+z,\tau)-w(y,\tau)|\omega_r(z)\, dy\, dz
   \\
\phantom{\lim_{\alpha_0\downarrow 0} R_{w,x}=}
 {}\le |w|_{L^{\infty}(0,T;BV({\mathbb R}^d))}
   \int_{{\mathbb R}^d} |z|\omega_r(z)\,dz
   \le C_1 r,
\end{gather*}
where $C_1:= |w|_{L^{\infty}(0,T;BV({\mathbb R}^d))}$.
Finally we have
\[
\lim_{\alpha_0\downarrow 0} R_{w^{\varepsilon},w} =
  \int_{{\mathbb R}^d}|w^{\varepsilon}(x,\tau)-w(x,\tau)|\, dx
   -  \int_{{\mathbb R}^d}|w^{\varepsilon}(x,\nu)-w(x,\nu)|\, dx.
\]

In summary from \eqref{total_entropy} we obtain
the following approximation inequality
\begin{gather}
      \int_{{\mathbb R}^d}|w^{\varepsilon}(x,\tau)-w(x,\tau)|\, dx \nonumber\\
\qquad {}      \le \int_{{\mathbb R}^d}|w^{\varepsilon}(x,\nu)-w(x,\nu)|\, dx
      + C_1 r  +\lim_{r_0,\alpha_0\downarrow 0}
      \Bigl( \overline{R}_{\mathrm{visc}}
      +  R_{\mathrm{conv}} \Bigr).
   \label{Approx_ineq}
\end{gather}

We start with the estimation of $\overline{R}_{\mathrm{visc}}$,
which can be done as follows:
\begin{gather}
      \overline{R}_{\mathrm{visc}} \le
      \varepsilon\sum_{i=1}^d\iiiint\limits_{Q_T\times Q_T}
      \bigl|\partial_{x_i}w^{\varepsilon}\bigr|\,
      \psi_{\alpha_0}(t)\bigr|\partial_{x_i}\omega_r(x-y)\bigr|
      \rho_{r_0}(t-s)\, dt\, dx\, ds\, dy
      \nonumber\\
\qquad {} \overset{\alpha_0\downarrow 0}{\longrightarrow}
      \sum_{i=1}^d \int_{\nu}^{\tau}\int_{{\mathbb R}^d}\int_{{\mathbb R}^d}
      \bigl|\partial_{x_i}w^{\varepsilon}\bigr|\,
      \bigr| \partial_{x_i}\omega_r(x-y)\bigr|\, dx\, dy\, dt
     \nonumber\\
\qquad {} \le \varepsilon K/r\sum_{i=1}^d \int_{\nu}^{\tau}\int_{{\mathbb R}^d}
      \bigr|\partial_{x_i}w^{\varepsilon}\bigr|\, dt\, dx
      \leq  \varepsilon T K/r |w^{\varepsilon}|_{L^{\infty}(0,T;BV({\mathbb R}^d))}
      \le C_2 T\varepsilon /r,
   \label{R_eps}
\end{gather}
where $K:=\int_{{\mathbb R}^d}\bigl|\delta'(\sigma)\bigr|\,d\sigma$
and $C_2:=K|w^{\varepsilon}|_{L^{\infty}(0,T;BV({\mathbb R}^d))}$.

Before we continue with the estimation of
$R_{\mathrm{conv}}$ we write $I_{\mathrm{conv}}=
I_{\mathrm{conv}}^1 + I_{\mathrm{conv}}^2$, where
\begin{gather*}
   I_{\mathrm{conv}}^1 = \mathrm{sgn}(w^{\varepsilon}-w)
   \Bigl[\bigl(V(y)f(w)-V(x)f(w)\bigr)\cdot \nabla_x\phi
 \\
\phantom{I_{\mathrm{conv}}^1 =}{}    -\bigl(V(x)f(w^{\varepsilon})- V(y)f(w^{\varepsilon})\bigr)
   \cdot \nabla_y\phi\Bigr],   \\
I_{\mathrm{conv}}^2 = \mathrm{sgn}(w^{\varepsilon}-w)
\bigl(\mathrm{div}_y V(y)f(w^{\varepsilon})
   - \mathrm{div}_x V(x)f(w)\bigr)\phi,
\end{gather*}
so that
\begin{gather*}
R_{\mathrm{conv}}=R_{\mathrm{conv}}^1 +R_{\mathrm{conv}}^2, \\
R_{\mathrm{conv}}^1=
\iiiint\limits_{Q_T\times Q_T} I_{\mathrm{conv}}^1 \, dt\, dx\, ds\, dy,
\qquad R_{\mathrm{conv}}^2
= \iiiint\limits_{Q_T\times Q_T} I_{\mathrm{conv}}^2 \, dt\, dx\, ds\, dy.
\end{gather*}

We start by estimating
$R_{\mathrm{conv}}^1$.  To this end introduce
\[
F(w^{\varepsilon},w):=\mathrm{sgn}(w^{\varepsilon}-w)\bigl[f(w^{\varepsilon})-f(w)\bigr]
\]
and observe that since $\nabla_y \phi =-\nabla_x \phi$,
\[
R_{\mathrm{conv}}^1=\iiiint\limits_{Q_T\times Q_T}
\Bigl((V(x)-V(y))F(w^{\varepsilon},w)\Bigr)\cdot\nabla_x \phi
\, dt\, dx\, ds\, dy.
\]
The function $F(\cdot,\cdot)$
is locally Lipschitz continuous in both variables
and the common Lipschitz constant equals $\mathrm{Lip}(f)$.
Since $w^{\varepsilon}\in L^{\infty}(Q_T)\cap L^{\infty}(0,T;BV({\mathbb R}^d))$,
$\nabla_xF(w^{\varepsilon},w)$ is a finite measure and
\[
\iint\limits_{Q_T}\bigl|\partial_{x_i}F(w^{\varepsilon},w)\bigr|\, dt\, dx\le
\mathrm{Lip}(f)\iint\limits_{Q_T}
\bigl|\partial_{x_i}w^{\varepsilon}\bigr|\, dt\, dx, \qquad i=1,\dots,d.
\]
Integration by parts thus gives
\begin{gather*}
   R_{\mathrm{conv}}^1 = -\underbrace{\iiiint\limits_{Q_T\times Q_T}
   \Bigl(\mathrm{div}_x V(x) F(w^{\varepsilon},w)
   \psi_{\alpha_0}(t)\omega_r(x-y)\rho_{r_0}(t-s)
   \, dt\, dx\, ds\, dy}_{R_{\mathrm{conv}}^{1,1}}
   \\  \phantom{R_{\mathrm{conv}}^1 =}
   {}- \underbrace{\iiiint\limits_{Q_T\times Q_T}
   (V(x)-V(y))\cdot \nabla_x F(w^{\varepsilon},w)
   \psi_{\alpha_0}(t)\omega_r(x-y)\rho_{r_0}(t-s)
   \, dt\, dx\, ds\, dy}_{R_{\mathrm{conv}}^{1,2}}.
\end{gather*}
For $R_{\mathrm{conv}}^{1,2}$ we calculate as follows:
\begin{gather*}
   \!\!\bigl|R_{\mathrm{conv}}^{1,2}\big| \! \le \,
    \mathrm{Lip}(f)\sum_{i=1}^d\!\iiiint\limits_{Q_T\times Q_T}\!
   \bigl|V_i(x)-V_i(y)\bigr|\, \bigl|\partial_{x_i} w^{\varepsilon}\bigr|
   \psi_{\alpha_0}(t)\omega_r(x-y)\rho_{r_0}(t-s)\, dt\, dx\, ds\, dy\!
   \\ \qquad {}\overset{\alpha_0\downarrow 0}{\longrightarrow}
   \mathrm{Lip}(f)\sum_{i=1}^d\int_{\nu}^{\tau}\int_{{\mathbb R}^d}\int_{{\mathbb R}^d}
   \bigl|V_i(x)-V_i(y)\bigr| \,
   \bigl|\partial_{x_i} w^{\varepsilon}\bigr|\, \omega_r(x-y)\, dx\, dy\, dt
   \\ \qquad {}\overset{(z:=x-y)}{=}\mathrm{Lip}(f) \sum_{i=1}^d\int_{\nu}^{\tau}
   \int_{{\mathbb R}^d}\int_{{\mathbb R}^d}
   \bigl|V_i(y+z)-V_i(y)\bigr|\,
   \bigl|\partial_{y_i} w^{\varepsilon}(y+z,t)\bigr|\omega_r(z)\, dz\, dy\, dt
   \\ \qquad {} \le \mathrm{Lip}(V)\mathrm{Lip}(f)
   \sum_{i=1}^d\int_{\nu}^{\tau}\int_{{\mathbb R}^d}
    \int_{{\mathbb R}^d}|z|
   \bigl|\partial_{y_i} w^{\varepsilon}(y+z,t)\bigr|\omega_r(z)\, dz\, dy\, dt
   \\ \qquad {} \le T \mathrm{Lip}(V)\mathrm{Lip}(f)
   |w^{\varepsilon}|_{L^{\infty}(0,T;BV({\mathbb R}^d))}
  \int_{{\mathbb R}^d}|z|\omega_r(z) \, dz
   \\ \qquad \le r T \mathrm{Lip}(V)\mathrm{Lip}(f)
   |w^{\varepsilon}|_{L^{\infty}(0,T;BV({\mathbb R}^d))}
   \le C_3 T r,
\end{gather*}
where $\mathrm{Lip}(V):=\max\limits_{i=1,\dots,d}\mathrm{Lip}(V_i)$ and
$C_3:=\mathrm{Lip}(V)\mathrm{Lip}(f)
|w^{\varepsilon}|_{L^{\infty}(0,T;BV({\mathbb R}^d))}$. Note that
we have used the Lipschitz regularity of
the velocity field $V$ (see~\eqref{Coeff_cond})
to get the desired result.

Regarding the term $R_{\mathrm{conv}}^2$, we firstly rewrite it as
\begin{gather*}
   R_{\mathrm{conv}}^2   =
   \underbrace{\iiiint\limits_{Q_T\times Q_T}
   \mathrm{div}_x V(x) F(w^{\varepsilon},w)
   \psi_{\alpha_0}(t)\omega_r(x-y)\rho_{r_0}(t-s)
   \, dt\, dx\, ds\, dy}_{R_{\mathrm{conv}}^{2,1}}
   \\ \hspace*{-20pt}{} +\! \underbrace{\iiiint\limits_{Q_T\times Q_T}
   \mathrm{sgn}(w^{\varepsilon}-w)
\bigl(\mathrm{div}_y V(y)  - \mathrm{div}_x V(x)\bigr) f(w^{\varepsilon})
   \psi_{\alpha_0}(t)\omega_r(x-y)\rho_{r_0}(t-s)
   \, dt\, dx\, ds\, dy}_{R_{\mathrm{conv}}^{2,2}}.
\end{gather*}
We set
\[
D(x)=\mathrm{div} V(x),
\]
and keep in mind that $D\in BV\big({\mathbb R}^d\big)$ by
\eqref{Coeff_cond}. Now we estimate $R_{\mathrm{conv}}^{2,2}$ as
follows:
\begin{gather*}
   \bigl|R_{\mathrm{conv}}^{2,2}\bigr|
   \le \|f(w^{\varepsilon})\|_{L^{\infty}(Q_T)}
   \iiiint\limits_{Q_T\times Q_T}
   \bigl|D(y)  - D(x)\bigr|
   \psi_{\alpha_0}(t)\omega_r(x-y)\rho_{r_0}(t-s)\, dt\, dx\, ds\, dy
   \\ \phantom{\bigl|R_{\mathrm{conv}}^{2,2}\bigr|} {} \overset{\alpha_0\downarrow 0}{\longrightarrow}
   \|f(w^{\varepsilon})\|_{L^{\infty}(Q_T)}
   \int_{\nu}^{\tau}\int_{{\mathbb R}^d}\int_{{\mathbb R}^d}
    \bigl|D(y)  - D(x)\bigr|
   \omega_r(x-y)\, dx\, dy\, dt
   \\ \phantom{\bigl|R_{\mathrm{conv}}^{2,2}\bigr|} {} \overset{(z=x-y)}{=}
   \|f(w^{\varepsilon})\|_{L^{\infty}(Q_T)}
   \int_{\nu}^{\tau}\int_{{\mathbb R}^d}\int_{{\mathbb R}^d}
   \bigl|D(y)  - D(y+z)\bigr|
   \omega_r(z)\, dz\, dy\, dt
   \\ \phantom{\bigl|R_{\mathrm{conv}}^{2,2}\bigr|}
 \le T\|f(w^{\varepsilon})\|_{L^{\infty}(Q_T)}
   \bigl|D\bigr|_{BV({\mathbb R}^d)}
   \int_{{\mathbb R}^d} |z|\omega_r(z)\,dz
   \le C_4 T r,
\end{gather*}
where $C_4:=\|f(w^{\varepsilon})\|_{L^{\infty}(Q_T)}
   \bigl|D\bigr|_{BV({\mathbb R}^d)}$. Note that
we have used the $BV$ regularity of $\mathrm{div} V$
to get the desired result.
Since $R_{\mathrm{conv}}^{1,1}
= R_{\mathrm{conv}}^{2,1}$, we have
\begin{equation}
   \label{R_conv}
   R_{\mathrm{conv}}=R_{\mathrm{conv}}^1
   + R_{\mathrm{conv}}^2 \le C_5 T r, \qquad
   C_5=\max(C_3,C_4).
\end{equation}

Set $C_6=\max(C_1,C_2,C_5)$. Then from
\eqref{Approx_ineq}, \eqref{R_eps} and \eqref{R_conv} we get
\begin{gather}
      \int_{{\mathbb R}^d}|w^{\varepsilon}(x,\tau)-w(x,\tau)|\, dx
      \le \int_{{\mathbb R}^d}|w^{\varepsilon}(x,\nu)-w(x,\nu)|\, dx
      + C_6\left((1+T)r + \frac{T\varepsilon}{r}\right)
       \nonumber\\   \qquad\qquad\qquad\qquad\qquad\quad
     {}
      \overset{\nu\downarrow 0}{\longrightarrow}
      C_6\left((1+T)r + \frac{\varepsilon}{r}\right).
\end{gather}
By choosing $r=\sqrt{T\varepsilon}$ we immediately
obtain
\begin{equation}
   \label{L1_tmp}
   \int_{{\mathbb R}^d}|w^{\varepsilon}(x,\tau)-w(x,\tau)|\, dx
   \le C_7 \sqrt{T\varepsilon}
\end{equation}
for some constant $C_7$ independent of $\varepsilon$.
To obtain \eqref{Error_estimate}, we simply
integrate \eqref{L1_tmp} over $\tau\in (0,T)$.

\subsubsection*{Added in process}
After the main result of this paper was obtained,
we became aware of a paper by Eymard, Gallouet and
Herbin~\cite{Eymard_etal_I:2000} which
also proves an error estimate for viscous approximate
solutions. They, however, deal with a certain
boundary value problem with a divergence free velocity field
and obtain an error estimate of order~$\varepsilon^{\frac15}$.
As is the case herein, the proof in~\cite{Eymard_etal_I:2000}
does not rely on a continuous dependence estimate.

\subsection*{Acknowledgments}

This work was done
while the first author (Evje) was visiting
the Industrial Mathematics Institute at the University of
South Carolina. Part of this work was completed
while the second author (Karlsen) was visiting the
Department of Mathematics and
the Institute for Pure and Applied Mathematics (IPAM)
at the University of California, Los Angeles (UCLA).

\label{evije-lastpage}

\end{document}